\documentclass{article}
\usepackage{spconf,amsmath,graphicx, amssymb}
\usepackage{hyperref}
\usepackage{enumitem}
\setlist{nosep, leftmargin=14pt}

\usepackage{siunitx}

\title{Computational Reproducibility in Metabolite Quantification Applied to Short Echo Time {\it{in vivo}} MR Spectroscopy }

\name{G. Vila\,$^{1}$, A. Bonnet\,$^{1}$, F. Chauveau\,$^{2}$, S. Gaillard\,$^{1}$, F. Durand-Dubief\,$^{1,3}$, S. Pop\,$^{1}$, H. Ratiney\,$^{1}$}
\address{$^{1}$Univ Lyon, INSA‐Lyon, Université Claude Bernard Lyon 1, UJM-Saint Etienne, CNRS, Inserm,
\\CREATIS UMR 5220, U1294, Lyon, France. \\$^{2}$Univ. Lyon 1, Lyon Neuroscience Research Center, CNRS UMR5292, INSERM U1028, Lyon, France.
\\$^{3}$ Service de Neurologie, Hôpital Pierre Wertheimer, GHE, Hospices Civils de Lyon, F-69677 Bron, France.%
}

\begin{document}
%
\maketitle
\begin{abstract}

\textit{In vivo} metabolite quantification by short echo time MR spectroscopy is a challenge for which various methods have been proposed. In this study, the reproducibility of quantification outcomes is questioned at three distinct levels: (i) between-software (LCModel and cQUEST), (ii) within-software (with different parameter sets), and (iii) across software executions (when the fitting algorithm uses random seeds, like cQUEST). After running multiple quantification tasks on a dedicated platform (VIP), metrics from Bland-Altman analysis were used to assess the variability of outcomes in signals acquired on a lysolecithin rat model, from a study on demyelination. Results show substantial variations at the three levels, allowing for more potent analyses than from a single parameter set / single software point of view.

\end{abstract}
\section{Introduction}
\label{sec:intro}

\subsection{Background}
\label{ssec:background}
\textit{In vivo}  MR spectroscopy is currently considered the only tool allowing for non-invasive measurement of biochemical compounds (\textit{e.g}, metabolites, lipids). 
Since pioneering publications presenting \textit{in vivo} spectra of rat brain \cite{pfeufferVivoNeurochemicalProfile1999},
signals acquired at short echo time in the brain have been of great interest because a large number of metabolites contribute to those signals. The main accessible metabolites are: N-acetylaspartate (NAA), creatine/phospho creatine (Cr+PCr), cholines (PCho+GPC), glutamate (Glu), glutamine, $\gamma$-amino-butyric acid (GABA) and taurine (Tau).

 MR spectroscopy quantification consists in analyzing the acquired signal to estimate metabolite concentration. In this process, the main challenges lie in the high number of metabolites with overlapping spectral components, the presence of a background of macromolecules, and the weak signal/noise ratios (since metabolites are present at concentrations \num{e4} to \num{e5} times smaller than water). The most widely developed quantification methods are based on parametric adjustment, using classical mathematical optimization methods like Least Squares. 
The \textit{in vivo} MR spectroscopy community has proposed several methods for brain metabolite quantification (see \href{https://mrshub.org/software/}{https://mrshub.org/software/}). 

\subsection{Reproducibility crisis}
\label{ssec:reproducibility}


 However, such methodological diversity challenges one's ability to reproduce numerical results across software. This statement is consistent with a growing awareness in the scientific community of the reproducibility issues associated with many studies \cite{baker500ScientistsLift2016}, including neuroimaging \cite{poldrackScanningHorizonTransparent2017, botvinik-nezerVariabilityAnalysisSingle2020}. This paper focuses on \emph{computational} reproducibility, which can be defined as one's ability to obtain identical results (or at least, equivalent conclusions), by applying identical analyses to the same set of experimental data.   
 
 In metabolite quantification, reproducibility issues can been observed at three levels. (i) High flexibility in the computation method ({\it e.g.}, model and parameters) leads to high variability in the final results \cite{marjanskaResultsInterpretationFitting2022}. (ii) Beyond the chosen method, the analysis software leaves its own footprint \cite{bhogal1HMRSProcessing2017}. (iii) Within the same software, the non-deterministic parts of the fitting process lead to impactful variability in the numerical outcomes. Since there is no ground truth for {\it in vivo} neurochemical profiles, these sources of variability should be controlled when carrying out a MR spectroscopy study.

\subsection{Motivations and outline}
\label{ssec:outline}
To our knowledge, however, no study to date has simultaneously controlled for these three sources of variability. 
This paper aims to capture them on the same MR spectroscopy dataset, by comparing the quantification outcomes:
\begin{itemize}
    \item between \textbf{two quantification algorithms}, accounting for inter-software variability;
    \item between \textbf{two sets of parameters} for each software, accounting for model flexibility; 
    \item across \textbf{multiple executions} for each software-parameter set, accounting for random seeds in the fitting process.
\end{itemize}

\section{Materials and Methods}
\label{sec:materials}

\subsection{Experimental data}
\label{ssec:data}

The signals used in this study were acquired on horizontal 11.7T Bruker Biospec. The STEAM sequence was used with the following parameters: TE 2ms, TM 10 ms, TR 3.5 s, spectral width 5464Hz, number of acquired complex points 2048, number of accumulations 256 resulting in 15 min of acquisition time, VOI=$12{\mu}$L. The $^1$H NMR spectra were acquired from two distinct VOIs located on the corpus callosum of brains of Lysolecithin rat models. This animal model consists in an intracerebral injection of lysophosphatidylcholine \cite{zhang2019evaluation} which produces focal demyelination (lesion). With help from T2-weighted anatomical imaging, the first voxel (Vox1) was located on the supposed lesion; and the second one (Vox2) was placed in the contralateral hemisphere, supposedly preserved. From this database, a total of 32 signals (\textit{i.e.,} 16 from voxel Vox1 and 16 from voxel Vox2) were submitted to our quantification experiments. 
All experiments were performed according to procedures approved by the local ethic committee for animal experimentation. 

\subsection{Quantification models}
\label{ssec:models}

\subsubsection{Software}
\label{ssec:quantificationSoftware}

Two quantification methods were used in this study: LCModel \cite{provencherEstimationMetaboliteConcentrations1993} and cQUEST\cite{ratiney2010semi}. LCModel is one of the most widely used method for the analysis of short-echo time 1H MRS data. It performs a model fitting in the frequency domain using, for the model function, a linear combination of metabolites signals, to which is added lipid and macromolecules contributions. The whole process is regularized using splines to account for different baseline and lineshapes. 

cQUEST is an alternative approach, with similar model description (sum of weighted metabolite signals, macromolecules and lipid contribution) but performing the parameter fitting in the time domain, and allowing or not the handling of the baseline with HLSVD \cite{pijnappel1992svd}. Compared to LCModel, no regularization term is employed. In the implementation of cQUEST used in this paper, a multiple starting values strategy has been implemented with a random seed.

\subsubsection{Parameter sets}
\label{ssec:quantificationParameters}
For both LCModel and cQUEST, a set of parameters must be provided by the user to define the fit options, prior knowledge and the constraints on the parameters to fit.  All signals have been preprocessed, phased and normalized with the unsuppressed water peak amplitude beforehand. For both methods, the same metabolite basis set was used and the signal from the macromolecules (MM) was voluntarily parameterized in the same way with Gaussian lines (modeled from one MM acquired with metabolite nulling) and constraints on the amplitudes of these peaks were set. LCModel necessarily adds splines for the baseline to correct for inaccuracies from approximative MM models and regularizes its fit. 

\begin{figure}[htb]
\begin{minipage}[b]{1.0\linewidth}
  \centering
  \centerline{\includegraphics[width=8.8cm]{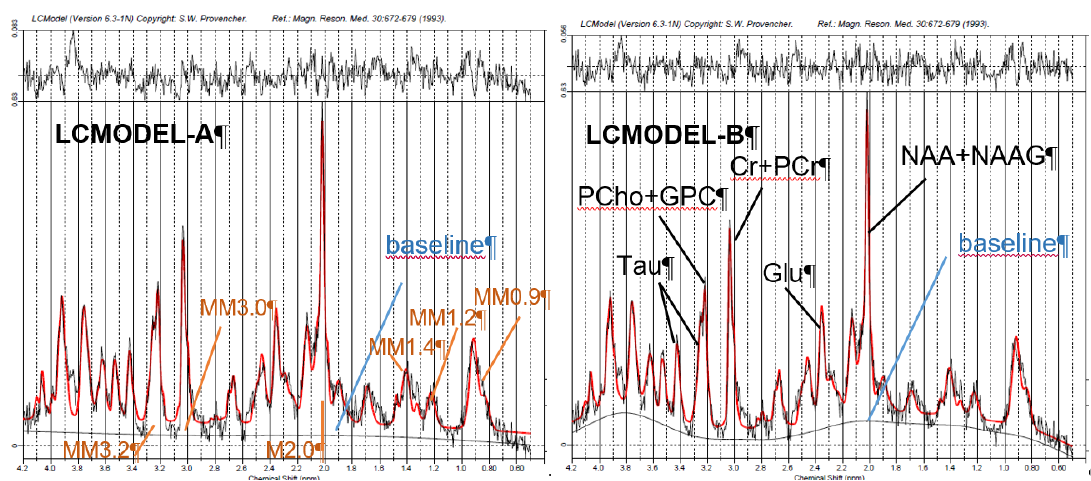}}
\end{minipage}
\vspace{-0.5cm}
\caption{LCModel results, with in red the fitted model overlaid on the raw spectrum in black for parameter set A and B. Annotations indicate baseline  macromolecules (MM) peak positions, and  metabolites.}
\label{fig:exampleParameterSets}
\end{figure}
Two solutions were computed for both methods, by playing on parameters acting on the baseline accounting (see Figure \ref{fig:exampleParameterSets}): no baseline (cQUEST-A) or baseline as flat as possible (LCModel-A) vs flexible baseline, HLSVD-modelled for cQUEST-B, and default option baseline modeling for LCModel-B.
\subsection{Reproducibility experiment}
\label{ssec:experiment}

\textbf{Experimental design $-$}
Mirroring the reproducibility questions raised in $\S$\ref{ssec:outline}, each of the 32 signals was submitted to metabolite quantification : (i) by both software (LCModel \& cQUEST), (ii) with both of their parameter sets (A \& B), on (iii) $n_{exec}$=30 distinct executions. This way, each metabolite was quantified 2*2*30=120 times for each signal in the dataset. Executions during which the quantification algorithm did not reach convergence were discarded, using indicators from both LCModel and cQUEST. This process left between 26 and 30 executions per signal. 

\textbf{Execution Environment $-$}
To perform this reproducibility experiment, both LCModel and cQUEST were containerized and made available on the Vitual Imaging Platform (VIP, see $\S$ \ref{sec:acknowledgments}). While in production conditions VIP uses distributed and heterogeneous resources, these experiments were run under controlled conditions to ensure that no variability was introduced in the outputs by the computing infrastructure or OS. We used two virtual machines, deployed on the same server and with the same OS version (CentOS 7 kernel $3.10.0-1160.76.1.el7.x86\_64$). 

\subsection{Measures of agreement}
\label{ssec:measuresOfAgreement}

A typical quantification result is the concentration $x(m,s,q,e)$ of metabolite $m$ in signal $s$, estimated by  method $q$ $\in$\{cquest-A, cquest-B, lcmodel-A, lcmodel-B\} during execution $e$. The purpose of this study is to measure the agreement in $x$ between two methods ($q$, $q'$); or across distinct executions $e$.

\subsubsection{Variability across software and parameter sets}
\label{ssec:method_blandAltman}
First, inter-execution variability was neutralized by using the mean concentration $x(m,s,q)$ over all executions $e$.
On this basis, a typical way to compare two measurements $q$ and $q'$ is Bland-Altman analysis \cite{martinblandSTATISTICALMETHODSASSESSING1986}, in which the differences $x_{m,s}(q)-x_{m,s}(q')$ are computed and displayed against the mean of both values. 

We call $bias(m,qq')$  the mean difference across signals, between methods $q$ and $q'$, for metabolite $m$; and $CI^{95}(m,qq')$ the 95\% confidence interval computed from the standard deviation of these differences. A simple way to assess if there is a strong disagreement between $q$ and $q'$ for metabolite $m$ is to compare both values in the following ratio:
\begin{equation}
Z^{95}(m,qq') = bias(m,qq') / CI^{95}(m,qq')
\label{eq:zscore}
\end{equation}
This modified Z-score was computed for each pair of quantification methods $qq'$ and each metabolite $m$.

\subsubsection{Variability across software executions}
\label{ssec:method_crossExec}

Inter-execution variability can be seen as the dispersion of $x$ around a "true" value $x_t$ : $x(e)=x_t+\epsilon(e)$, where $\epsilon(e)$ is the measurement error for execution $e$. A simple way to measure $\epsilon$ is to compute the residuals of a linear model, as in \cite{riemannAssessmentMeasurementPrecision2022}.

Metabolite concentration $x$ was assumed to depend on the signal $s$ ("$x \sim s$"); and software execution $e$ to produce a random effect. This model was fit separately on the signals from voxels Vox1 and Vox2.
For each metabolite $m$, method $q$ and voxel $v$, a root-mean-squared error (RMSE) was computed from the residuals $r(s,e)$ of the corresponding linear model:
\begin{equation}
RMSE(m,q,v) = 
    \sqrt{ 
        \underset{e, s \in v}{ \textrm{mean}}
        \left(
            r_{m,q}\left(s,e\right)^2
        \right)
        }
\label{eq:rmse}
\end{equation}

We used bootstrapping over signals to measure the influence of the dataset on the RMSE. The signal set used to fit the linear model was resampled randomly; and model fitting was run on this new dataset to give new residuals. This process was repeated 1000 times to give as many values of RMSE. 

This inter-execution measure of variability was compared with the traditionnal Cramér-Rao bound (CRB, the theoretical lower bound on the variance of an unbiased estimator), computed by both software for each metabolite concentration.

\subsubsection{Preservation of the findings}
\label{ssec:method_preservationFindings}

Once uncertainty measures have been performed, a final question to be answered is whether such variability may impact the results of a biological investigation on the dataset. 

As stated in \S \ref{ssec:data}, our signals should show variations in metabolite concentrations $x_{m,q}$ between the two voxels (Vox1 \& Vox2). For each execution $e$, a Wilcoxon sign-rank test was run on quantification outputs $x_{m,q,e}(s)$ between $s\in$ Vox1 and $s\in$ Vox2. In this process, we were not interested in the test results but in their consistency across all executions.

\section{Results}
\label{sec:results}

\subsection{Variability between parameter sets}
\label{ssec:betweenParameters}

\begin{figure}[htb]

\begin{minipage}[b]{.36\linewidth}
  \centering
  \centerline{\includegraphics[height=4.2cm]{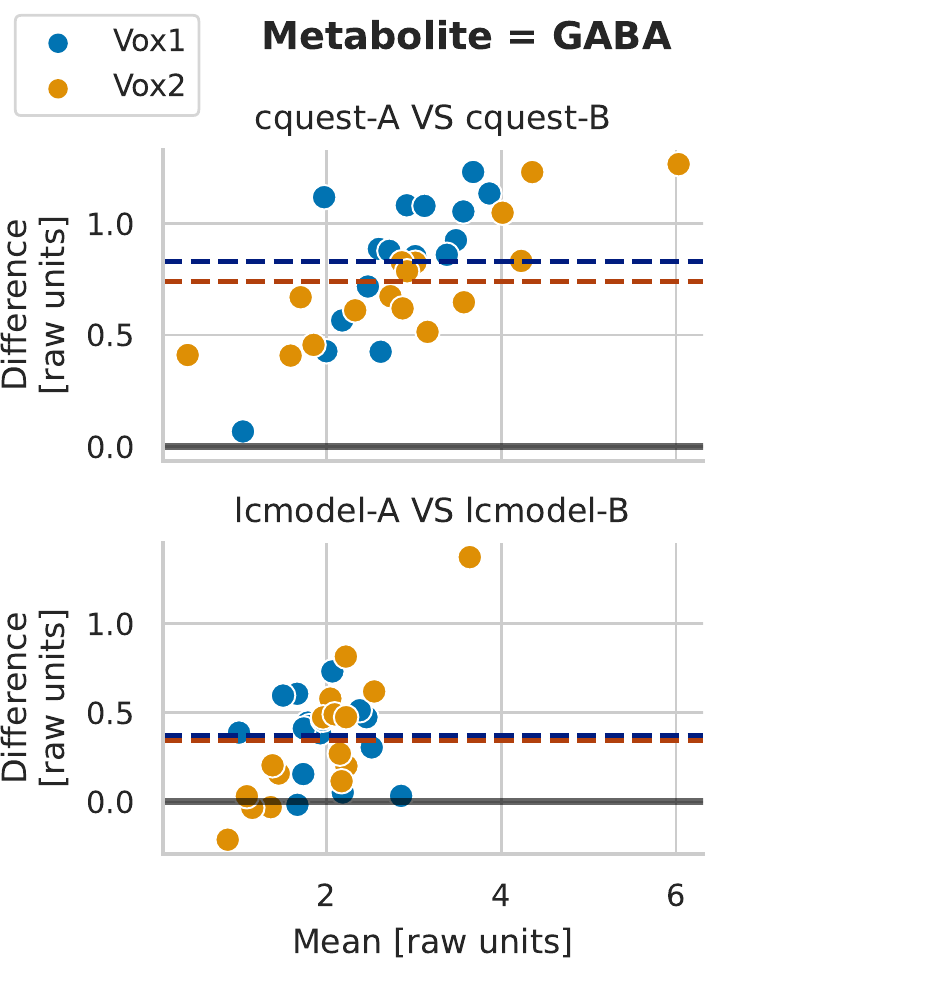}}
\end{minipage}
\begin{minipage}[b]{.30\linewidth}
  \centering
  \centerline{\includegraphics[height=4.2cm]{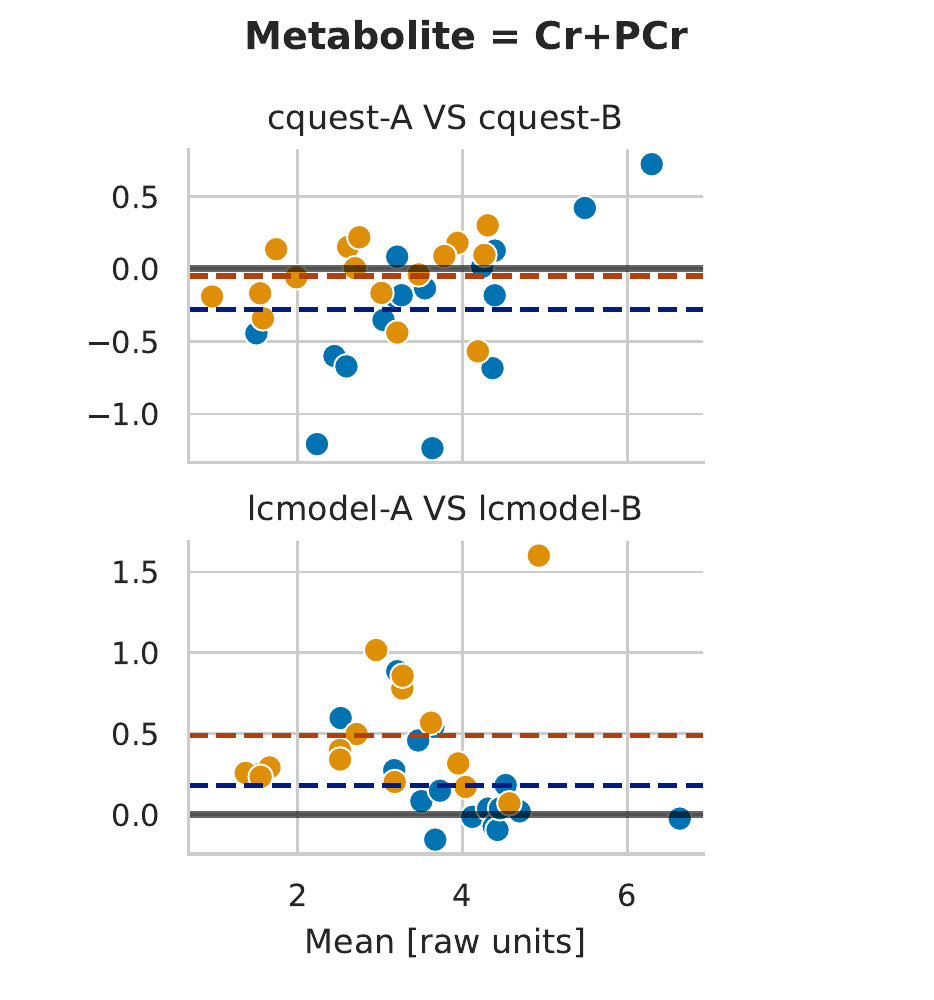}}
\end{minipage}
\begin{minipage}[b]{.30\linewidth}
  \centering
  \centerline{\includegraphics[height=4.2cm]{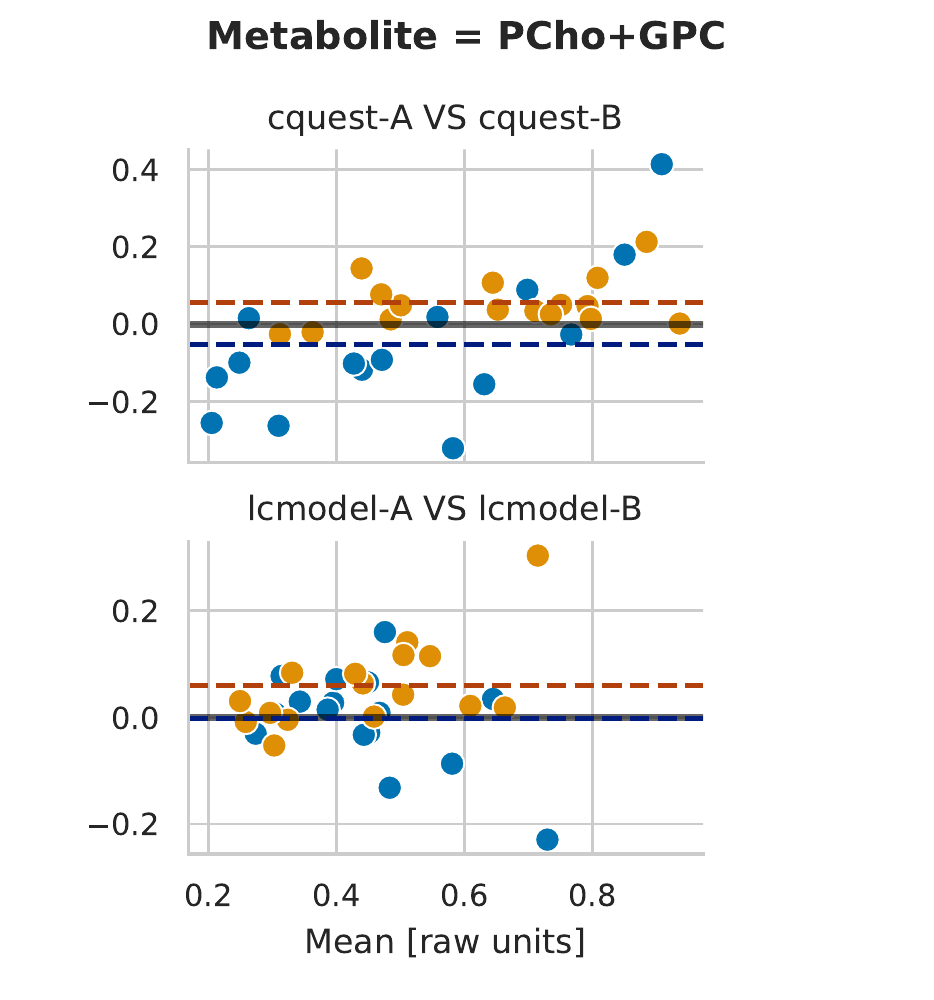}}
\end{minipage}

\vspace{-0.5cm}
\caption{Bland-Altman plots showing between-parameter variations in the quantification results of the same sotfware (LCModel, cQUEST) for three metabolites.}
\label{fig:betweenParameters}
\end{figure}

On Bland-Altman plots of Figure \ref{fig:betweenParameters}, the between-parameter difference (y-axis) is displayed against the mean (x-axis). For each metabolite, the upper plot stands for cQUEST and the lower plot for LCModel. 
For each voxel, the bias (see \S \ref{ssec:method_blandAltman}) is depicted by a dash-dotted line (blue for Vox1, orange for Vox2) and can be compared to Zero (plain dark grey line).


Both LCModel and cQUEST show systematic disagreements between their parameter sets. This is obvious on the left-hand charts of Figure \ref{fig:betweenParameters} (metabolite GABA), where most of the points lie above the Zero line. When focusing on Vox2, this is also true for Cr+PCr within LCModel; and for PCho+GPC within both models. 

For the three metabolites, Vox2 interestingly shows higher biases but lower dispersion than Vox1. This reflects the expected higher biological variability of metabolite concentration in the damaged tissue (\textit{i.e.}, Vox1). The bias / dispersion ratios are displayed on first and last lines of Figure \ref{fig:biasHeatmaps}, where within-software comparisons are highlighted in bold font.

\subsection{Variability between quantification methods}
\label{ssec:betweenMethods}

Figure \ref{fig:biasHeatmaps} extends the previous analyses to six metabolites and the six possible combinations of quantification methods ($qq'$). 
Each matrix displays the modified Z-scores introduced in Equation \ref{eq:zscore} for a single voxel. In a classical Z-test, a value above 100\% would mean significant variations (this is not the focus of the present study, as sample sizes are small).

On the upper 3 lines of each matrix, method cquest-A (which does not use any baseline modeling) shows strong deviations from the other methods regarding the quantification outputs. On the two right-hand columns, one may also notice a divergence in the results between both voxels, especially for metabolite pair PCho+GPC (absolute Z-scores are all below 30\% for Vox1 and above 36\% for Vox2). As stated earlier (\S \ref{ssec:betweenParameters}), the lower Z-scores for Vox1 may be explained by a higher dispersion in the between-method differences.

\begin{figure}[htb]
\begin{minipage}[b]{1.0\linewidth}
  \centering
  \centerline{\includegraphics[width=8.8cm]{figures/inter-model/BiasHeatmaps.pdf}}
\end{minipage}
\vspace{-0.8cm}
\caption{Bias between quantification methods, expressed as a percentage of the 95\% confidence intervals in Bland-Altman plots. C-A/B: cquest-A/B ; L-A/B: lcmodel-A/B .}
\label{fig:biasHeatmaps}
\end{figure}

\subsection{Variability across software executions}
\label{ssec:betweenExec}

\begin{figure}[htb]
\begin{minipage}[b]{1.0\linewidth}
  \centering
  \centerline{\includegraphics[width=8.5cm]{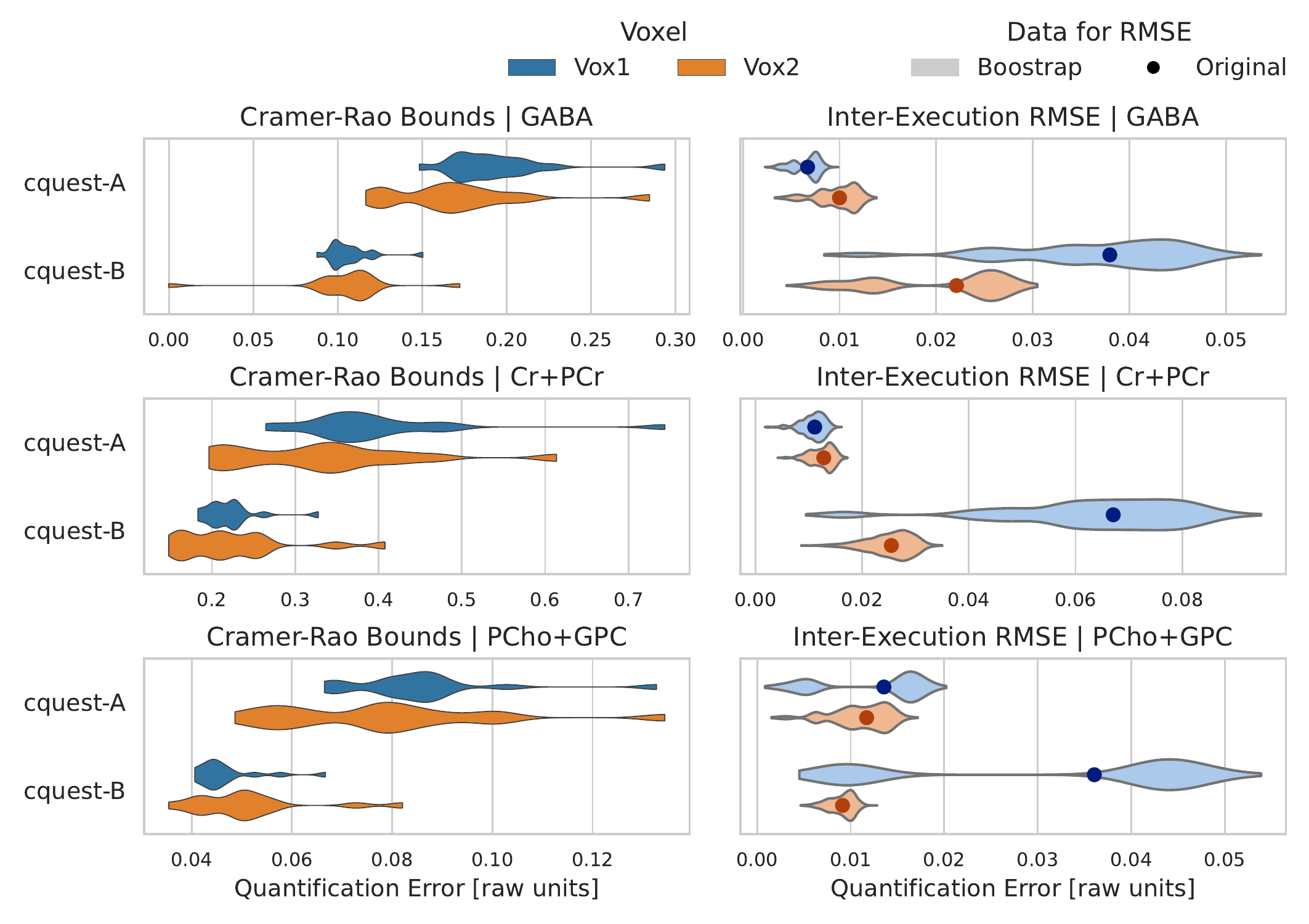}}
\end{minipage}
\vspace{-0.8cm}
\caption{Inter-execution variability compared to Cramér-Rao lower bounds for three metabolites.}
\label{fig:betweenExec}
\end{figure}

All previous analyses were performed on a single concentration value for each metabolite $m$, method $q$ and signal $s$.
To assess the quantification error on these values, Figure \ref{fig:betweenExec} displays the theoretical variance computed by the software (CRB, left-hand plots) next to our empirical measure of the inter-run variability (RMSE, right-hand plots).
In this study, only cQUEST could be explicitly designed to use random seeds; and LCModel outputs showed zero inter-run variability. Hence, the results are displayed for cQUEST only. 
\vspace{-0.2cm}

Two statements may be drawn from these results. First, the quantification errors are rather small here compared to the mean metabolite concentrations observed in Figure \ref{fig:betweenParameters}. Second, CRBs and RMSEs are inconsistent between each other: while CRBs (which are higher than RMSEs) decrease between cquest-A and cquest-B, RMSEs increase between both models. Unlike RMSEs, the CRBs do not seem to make a difference between Vox1 and Vox2. Therefore, CRB and the inter-run RMSE seem to provide two distinct measures of uncertainty on the quantification outcomes.

\subsection{Preservation of the findings}
\label{ssec:concdrawn}

\begin{figure}[htb]
\begin{minipage}[b]{1.0\linewidth}
  \centering
  \centerline{\includegraphics[width=6cm]{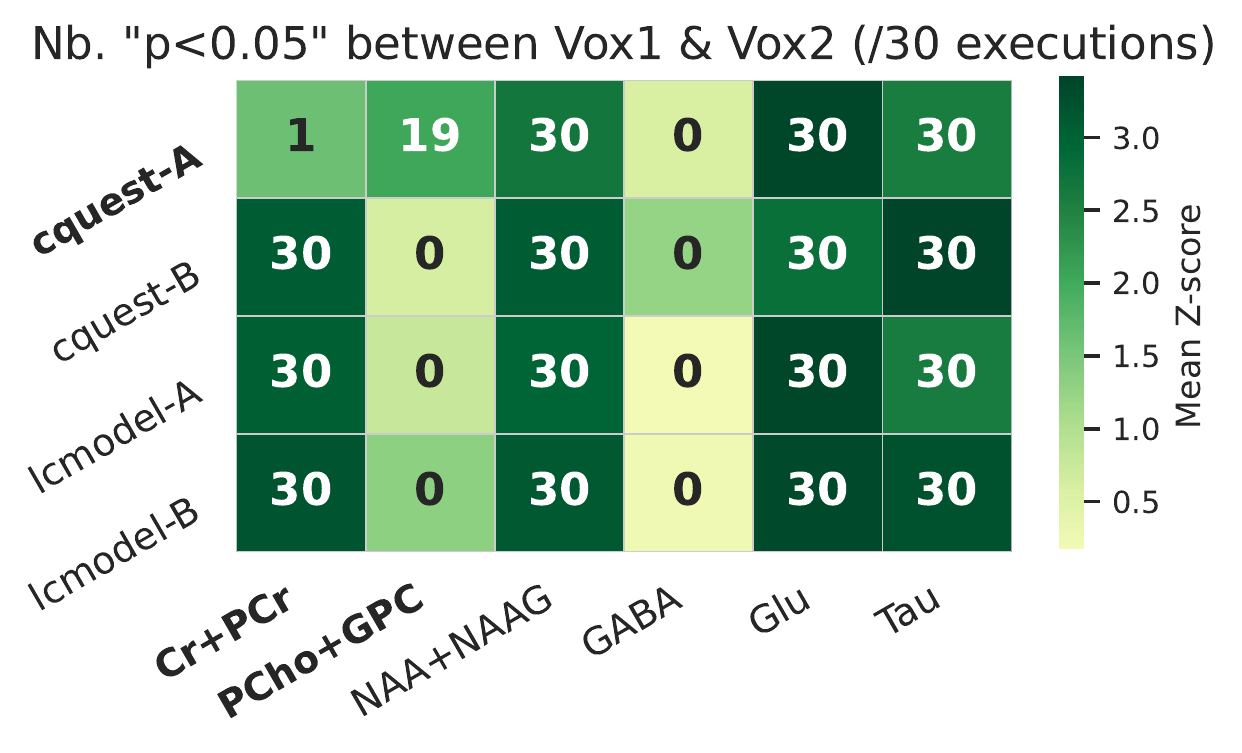}}
\end{minipage}
\vspace{-0.8cm}
\caption{Number of significant results (out of 30 quantification runs) obtained by the signed rank test, between Vox1 and Vox2. Colormap: average Z-statistic from the test.}
\label{fig:findingPreservation}
\end{figure}

Finally, the conclusions drawn from a signed rank test between Vox1 and Vox2 are summarized in Figure  \ref{fig:findingPreservation}. For each method and metabolite, the matrix displays the number of "significant" variations detected over the 30 quantification runs. Divergent results are highlighted in bold font. 

\vspace{-0.4cm}

Except for two metabolite pairs, the test delivered the same conclusions for each run regardless of the software or parameter set. For metabolites Cr+PCr and PCho+GPC, however, the test provided divergent results between cquest-A and the other methods; and divergent results between quantification runs.

\section{Discussion}
\label{sec:discussion}

Quantification of short echo time spectroscopy data is an ill-posed problem, for which the macromolecular contribution is not fully known; and the additional consideration of a baseline is very software dependent. Depending on the parameter set, the implementation of the inverse problem, or the regularization terms, multiple solutions can be obtained from the same dataset. This work proposes a methodology to monitor distinct levels of variability in the quantification outcomes, related to numerical implementations of the fitting method. 

Despite the strong biases observed between cQUEST and LCModel and within these models for some metabolites (\textit{e.g.} GABA), both models mostly agreed during hypothesis testing for differences between healthy and damaged brain regions. This was not the case, however, for the quantification of cholines and creatine (PCho+PCh and PCr+Cr). This may be caused by a strong interdependence between macromolecules, respectively at 3.2 ppm and 3.0 ppm and the peak of total Choline and total Creatine. That said, this diversity of solutions is \textit{in itself} a sign of something interesting to investigate. This fosters the use of multiple quantification methods, to overcome the limitations of a single software or parametrization in the fitting process. 

Finally, this reproducibility study was designed to be reproduced. The code and data used to deliver the present results are available on our GitLab repository: 

\href{https://gitlab.in2p3.fr/gael.vila/repro-spectro}{https://gitlab.in2p3.fr/gael.vila/repro-spectro}


\section{Compliance with ethical standards}
\label{sec:ethics}
All experiments were performed according to procedures approved by the local ethic committee for animal experimentation (\textit{Comité d’Expérimentation Animale de l’Université Claude Bernard Lyon 1 – CEEA-55}).

\section{Acknowledgments}
\label{sec:acknowledgments}

This work is part of the ReproVIP project\footnote{\url{https://anr.fr/Projet-ANR-21-CE45-0024}}, which aims at evaluating and improving the reproducibility of scientific results obtained with the Virtual Imaging Platform (VIP)\footnote{\url{https://vip.creatis.insa-lyon.fr}} in the field of medical imaging \cite{glatardVirtualImagingPlatform2013}. VIP is a web portal which enables researchers worldwide to easily execute medical imaging applications on a large scale computing infrastructure, by leveraging the computing and storage resources of the EGI federation \footnote{\url{https://www.egi.eu/egi-federation/}}. In October 2022, VIP counts 1445 registered users from 83 countries and approximately 20 applications available as a service.

This work was supported  by the French National Research Agency under Grant N°ANR-21-CE45-0024. It has also received funding from the Neurodis fondation. It was performed within the framework of LABEX PRIMES (ANR-11-LABX-0063) of Université de Lyon, within the program "Investissements d'Avenir" (ANR-11-IDEX-0007) and France Life Imaging (ANR-11-INBS-0006). This work is based on spectra acquisitions achieved in the PILoT facility, member of the infrastrucure France Life Imaging (ANR-11-INBS-0006).

\bibliographystyle{IEEEbib}
\bibliography{biblio}

\end{document}